\documentclass[12pt,twoside,a4paper]{amsart}
\usepackage{amssymb}
\date{\today}

\def\deg{\text{deg}\,}

\def\w{\wedge}

\def\dbar{\bar\partial}

\def\C{{\mathbb C}}
\def\P{{\mathbb P}}

\def\D{{\mathcal D}}

\def\lsigma{\text{\Large $\sigma$}}
\def\Hom{{\rm Hom\, }}
\def\codim{{\rm codim\,}}

\def\Im{{\rm Im\, }}

\def\Ker{{\rm Ker\,  }}

\def\E{{\mathcal E}}

\def\O{{\mathcal O}}

\def\L{{\mathcal L}}

\def\Re{{\rm Re\,  }}

\def\L{{\mathcal L}}

\def\be{\begin{equation}}
\def\ee{\end{equation}}

\newtheorem{thm}{Theorem}[section]
\newtheorem{lma}[thm]{Lemma}
\newtheorem{cor}[thm]{Corollary}
\newtheorem{prop}[thm]{Proposition}

\theoremstyle{definition}

\newtheorem{df}{Definition}

\theoremstyle{remark}

\newtheorem{preremark}{Remark}
\newtheorem{preex}{Example}

\newenvironment{remark}{\begin{preremark}}{\qed\end{preremark}}
\newenvironment{ex}{\begin{preex}}{\qed\end{preex}}

\numberwithin{equation}{section}

\title[Residue currents of  holomorphic 
\ldots ]
{Residue currents  of  holomorphic morphisms}

\begin{document}

\date{\today}

\author{Mats Andersson}

\address{Department of Mathematics\\Chalmers University of Technology and the University of G\"oteborg\\S-412 96 G\"OTEBORG\\SWEDEN}

\email{matsa@math.chalmers.se}

\subjclass{32 A 27, 32 H 02, 32 B 99}


\thanks{The author was
  partially supported by the Swedish Natural Science
  Research Council}

\begin{abstract}
Given a generically surjective holomorphic vector bundle morphism $f\colon E\to Q$,
$E$ and $Q$ Hermitian bundles,
 we construct a current
$R^f$ with values in $\Hom(Q,H)$, where $H$ is a certain  derived bundle,
and with support on the set $Z$ where $f$ is not surjective. The main property is that
if $\phi$ is a holomorphic section of $Q$, and $R^f\phi=0$, then
locally $f\psi=\phi$ has a holomorphic solution  $\psi$.
In the generic case also the converse holds. This gives a generalization
of the corresponding theorem for a complete intersection, due to
Dickenstein-Sessa and Passare.
We also present   results for polynomial mappings, related to 
M Noether's theorem and the effective Nullstellensatz. The construction of the current
is based on a generalization of the Koszul complex. By means of this complex one  can also
obtain new global estimates of solutions to $f\psi=\phi$, and as an example we give
new results related to the $H^p$-corona problem.
\end{abstract}



\maketitle

\section{Introduction}

Let $E$ and $Q$ be holomorphic Hermitian vector bundles of ranks $m$ and $r$,
respectively, over the $n$-dimensional complex manifold $X$, and let
$f\colon E\to Q$ be a generically  surjective holomorphic morphism.
Given a holomorphic section $\phi$ of $Q$ we are interested in holomorphic solutions $\psi$ to
$f \psi=\phi$. 
The basic results in this  area are the existence theorems due to Skoda 
in \cite{Sk1} and \cite{Sk2}, which are based on $L^2$-methods and 
complex geometry. They provide existence of global holomorphic solutions to
the equation $f\psi=\phi$ with $L^2$-estimates under appropriate geometric conditions
provided $f$ is pointwise surjective. However, applying these results to
$E$ restricted to $X\setminus Z$,  where
$$
Z=\{z;\  f(z)\ {\rm is\ not\  surjective}  \},
$$
also highly non-trivial local results at  $Z$ are obtained by these methods.

\smallskip

In this paper we  introduce a complex of bundles
$$
\cdots\to E_3\to E_2\to E\to Q\to 0,
$$
and  define a global  residue current 
$$
R^f=R^f_p+\cdots +R^f_\mu
$$
with support on $Z$,  $p=\codim Z$ and $\mu=\min(n,m-r+1)$,
where  $R^f_k$ is a $(0,k)$-current with values in $\Hom(Q,E_k)$.
It is not hard to see
(e.g., by using Gauss elimination) that  $p\le m-r+1$ with equality in the generic case; in this case,
thus $R=R_{m-r+1}$.
Our first   result concerns existence of local holomorphic solutions of
$f\psi=\phi$.

\begin{thm}\label{basal}
Let $E$ and $Q$ be holomorphic Hermitian vector bundles over a complex manifold $X$,  let
 $f\colon E\to Q$ be a holomorphic generically surjective morphism,  and let
 $R^f$ be the corresponding  residue current.
If $\phi$ is a holomorphic section of $Q$ such that $R^f\phi=0$, then
locally $f\psi=\phi$ has a holomorphic solution $\psi$.
\end{thm}

We have the following partial converse.

\begin{thm}\label{omv}
If $p=m-r+1$ and $f\psi=\phi$ has a holomorphic solution, then $R^f\phi=0$.
\end{thm}

If $p=m-r+1$ thus $f\psi=\phi$ has holomorphic solutions if and only if 
$R^f\phi=0$.
If  $r=1$ and $p=m$, it turns out that  if   $f=\sum f_j e_j$ in a local holomorphic frame  $e_j$,
then $R=R_m$ is equal to the classical Coleff-Herrera current
$$
T=\Big[\dbar\frac{1}{f_1}\w\ldots\dbar\frac{1}{f_m}\Big]
$$
times a non-vanishing section of $(\det E)\otimes Q^*$, see \cite{A2}. In this case 
therefore we  get back the  Dickenstein-Sessa-Passare theorem, \cite{DS} and \cite{P1},
stating  that $\phi$ belongs to the ideal $(f_1,\ldots,f_m)$ if and only if $\phi T=0$.


\smallskip

Instead of the usual norm $|\phi|$ of a section $\phi$ of $Q$ it is natural,
e.g., in view of the results in \cite{Sk2},  to introduce  the   stronger pointwise norm
$$
\|\phi\|^2=\det (ff^*)|f^*(ff^*)^{-1}\phi|^2=\langle \widetilde{ff^*}\phi,\phi\rangle,
$$
where $\widetilde{ff^*}=\det(ff^*)(ff^*)^{-1}$ is the smooth endomorphism on $Q$  whose matrix
 is the transpose of the comatrix of $ff^*$. 
By analyzing the singularity of $R^f$ we obtain the following 
sufficient size condition on  $\phi$ for annihilating the residue.

\begin{prop}\label{bsgen}
Let $f\colon E\to Q$ be a holomorphic generically surjective morphism. If $\phi$ is a holomorphic section of
$Q$ such that 
\begin{equation}\label{bsvillkor}
\|\phi\|^2\le C \det (ff^*)^{\min(n,m-r+1)},
\end{equation}
then  $R^f\phi=0$.
\end{prop}

As an immediate consequence  we get the following 
generalization of  the Brian\c con-Skoda theorem.

\begin{thm}\label{genbri}
If  $\phi$ is a holomorphic section of
$Q$ such that \eqref{bsvillkor} holds, then locally $f\psi=\phi$ has a holomorphic solution  $\psi$.
\end{thm}

In the case when $r=1$, \eqref{bsvillkor} means precisely   that  $|\phi|\le C |f|^{\min(n,m)}$,
 and the conclusion is then that
$\phi$ is locally in the ideal $(f)$ generated by $f$
(i.e., the ideal generated by $f_j$ if $f=\sum f_j e_j$ in some local holomorphic frame $e_j$). 
This immediately implies the classical
Brian\c con-Skoda theorem, \cite{BS}, which  states that $\phi^{\min(m,n)}$ belongs to $(f)$
if $|\phi|\le C|f|$. 
For $m-r+1\le n$,  Theorem~\ref{bsgen} also follows directly
from Skoda's $L^2$-estimate, \cite{Sk2},  
but when $m-r+1>n$,  the $L^2$-estimate only gives the conclusion if the power 
of  the right hand side of \eqref{bsvillkor} is
$n+1$.  By an additional argument in  the case when $m>n$,  the classical Brian\c con-Skoda theorem  follows from  
the $L^2$-estimate;  the case when $m>n$ and $r>1$ can be reduced to the classical result by means
of the Fuhrmann trick, see Section~\ref{hpcorona}.

\smallskip
Demailly has extended Skoda's $L^2$-theorems to $\dbar$-closed sections,
see \cite{Dem0} and \cite{Dem}. Our method
also admits such an extension of the local result.

\begin{thm}\label{basal2}
Assume that $\phi$ is a smooth $\dbar$-closed $(0,q)$-section of $Q$. If 
$R^f\phi=0$, then
locally $f\psi=\phi$ has   $\dbar$-closed (current)
solutions  $\psi$.
\end{thm}

For degree reasons we see that $R^f\phi=0$ if $q>n-p$. Thus we get

\begin{cor}
If $\phi$ is any  smooth $\dbar$-closed $(0,q)$-form with values in $Q$, and
$q>n-p=\dim Z$, then  locally $f\psi=\phi$ has  $\dbar$-closed current solutions.
\end{cor}

In analogy with Theorems~\ref{bsgen} and \eqref{genbri}  we also have

\begin{thm}\label{bsgen2}
If $\phi$ is a smooth $\dbar$-closed $(0,q)$-form with values in $Q$ such that 
\begin{equation}\label{bsvillkor2}
\|\phi\|^2\le C\det(ff^*)^{\min(n-q,m-1+r)},
\end{equation}
then $R^f\phi=0$,
and locally there are integrable  $\dbar$-closed solutions $\psi$ to $f\psi=\phi$.
\end{thm}

\smallskip

We can also obtain global  results and first we 
turn our attention to polynomial ideals and  generalize the
approach in \cite{A3}.  
Let $[z]=[z_0,\ldots,z_n]$ be  homogeneous coordinates on $\P^n$, and let 
$z'=(z_1,z_2,\ldots,z_n)$ be the standard coordinates in the standard affinization
$\C^n\simeq\{[z], z_0\neq 0\}$.
Let $P$ be a polynomial mapping $\C^n\to \Hom(\C^m,\C^r)$ with columns $P^j$ such that 
$\deg P^j \le d_j$, $j=1,\ldots,m$.
If $f$ is the matrix whose columns are the  $d_j$-homogenized forms 
$f^j(z)=z_0^{d^j}P^j(z'/z_0)$ in $\C^{n+1}$,  then $f$  defines a morphism
$$
f\colon \bigoplus_1^m\O(-d_j)\to \C^r.
$$
Let $Z$ be the algebraic variety in $\P^n$ where $f$ is not surjective,
and let $R^f$ be the associated residue current with respect to the natural metric.

\begin{thm}\label{pellenom}
Assume that $P$ is a polynomial mapping as above, and 
let $\Phi$ be  a $r$-column of polynomials of degrees $\le \rho$.
Moreover, assume that
\begin{equation}\label{villkor}
m\le n+r-1 \quad {\rm or}\quad \rho\ge \sum_{j=1}^{n+r}d_j -n,
\end{equation}
where $d_1\ge d_2\ge \ldots \ge d_m$.  
If
$
R^f\phi=0,
$
then there are polynomials $Q_j$ such that 
$\sum_1^m  P^jQ_j=\Phi$,
and
$\deg P^jQ_j\le \rho$. 
\end{thm}


\smallskip
\begin{cor}
Assume that $Z$ is empty. Then we can find a matrix
$Q$ of polynomials with rows $Q_k$ 
such that  $PQ=\sum P^kQ_k=I_r$, and
$\deg P^kQ_k\le   \sum_{j=1}^{n+r}d_j-n$.
\end{cor}

This is a generalization of a classical theorem of Macaulay,
\cite{Macaul}.

\smallskip
\begin{cor}\label{groda}
Let $\Phi$ be  a column of polynomials, $\deg\Phi\le\rho$, and  let $\phi$ be  its
$\rho$-homogenization. If 
$$
\|\phi\|^2\le C\det(ff^*)^{\min(n,m-r+1)}
$$
in $\P^n$, and \eqref{villkor} is fulfilled, then
$PQ=\Phi$ has a solution with
$\deg P^kQ_k\le \rho$.
\end{cor}

\smallskip

Assume that $P$ is pointwise surjective in $\C^n$,
and let $P^j=(P_1^j,\ldots,P_r^j)^t$.
By the local Lojasiewicz inequality there is a constant  $M$ such that 
\begin{equation}\label{olik}
\sum_{|I|=r}'\frac{|\det(P^{I_j}_k(z'))|^2}{(1+|z'|^2)^{\sum_1^rd_{I_j}}}
\ge C\frac{1}{(1+|z'|^2)^M},
\end{equation}
where the sum is over increasing multiindices. 

\begin{cor}
Assume that  $P\colon\C^n\to\Hom(\C^m,\C^r)$ is surjective in $\C^n$, $\deg P\le d$
and that \eqref{olik} holds. 
Then there is a matrix $Q$ of polynomials  such that 
$ PQ=I_r$
and 
$\deg P^jQ_j\le M\min(n,m-r+1)$.
\end{cor}

From Kollar's famous theorem, \cite{Koll}, we can get an  estimate of $M$.
For simplicity we assume $d_j=d$ for all $j$.

\begin{prop}\label{kollest}
If $d_j=d$ for all $j$, then the inequality \eqref{olik} holds with
$$
M=(rd)^{\min(n,m!/(m-r)!r!)},
$$
provided that $rd\ge 3$.
\end{prop}

It should be pointed out that the bound 
$$
\deg Q+d\le \min(n,m-r+1)M
$$  
we obtain in this way for a solution to $PQ=I_r$ (even)
is not optimal when $r=1$. It is proved in \cite{Koll} that one actually have
$\deg Q+d\le M$ when $r=1$. We do not know if it is  possible to modify Kollar's proof
as to include the case $r>1$ directly and  get a sharper bound.

\smallskip
Now assume that $p=m-r+1=n$ and that $Z$ is contained in $\C^n$; thus a discrete set.
Moreover, assume that $\Phi=P Q$ is solvable in $\C^n$. Then it follows from
Theorem~\ref{omv} that $ R^f\phi=0$ in $\C^n$, and hence $R^f\phi=0$ in $\P^n$, and since
\eqref{villkor}  is fulfilled we can take $\rho=\deg \Phi$. Therefore there is a solution to
$PQ=\Phi$ such that  $\deg P^jQ_j\le \deg \Phi$. 
When $r=1$ this is a classical  theorem due to  Max Noether, \cite{Noe}.
We have the following generalization that  appeared in \cite{A3} in the case
$r=1$;  however  we suspect that this case  could be proved algebraically,  e.g.,   by the methods
in \cite{Hick}.

\begin{thm}\label{noethergen}
Assume that $P\colon\C^n\to\Hom(\C^m,\C^r)$ and that $p=m-r+1$ and that
$Z$ has no irreducible component contained in the hyperplane at
infinity.
Moreover, assume that $\Phi=P Q$ is solvable in $\C^n$. Then there 
is a solution $Q$ such that
$\deg P^jQ_j\le \deg \Phi.$
\end{thm}

\smallskip
We can  also obtain new global results in open bounded domains even when $f$ is pointwise surjective,
and as an example we present in Section~\ref{hpcorona} new sharpened estimates of solutions to the
$H^p$-corona problem in a strictly pseudoconvex domain.

\tableofcontents

\section{A generalized Koszul complex}\label{nykomplex}

Let $f\colon E\to Q$ be  a   holomorphic morphism as above. Assume that we have a complex
\begin{equation}\label{complex}
\cdots \to E_3\to E_2\to E\to Q\to 0.
\end{equation}
of vector bundles where all the morphisms, which we denote by $\delta$, are holomorphic.
Let $E_0=Q$, $E_1=E$, and let 
$$
H=\bigoplus_{k=0}^\infty E_k.
$$
We will also consider  $(0,*)$-form-valued sections of $H$, i.e.,
sections of $T^*_{0,1}(X)\otimes H$.
We denote this space of sections by  $\E_{0,*}(X,H)$. Notice that
it is a module over the ring (algebra) $\E_{0,*}(X)$.
We extend the action of $\delta$ to sections of
$T^*_{0,1}(X)\otimes H$ by requiring that
\begin{equation}\label{extend}
\delta \xi\otimes w =(-1)^{\deg \xi} \xi\otimes \delta w
\end{equation}
if $\xi$ is a differential form.  Then  $\delta \dbar=-\dbar\delta$.
\smallskip

Now suppose that we have $(0,k-1)$-forms, or currents, $v_k$ with values in $E_k$, $k\ge 1$,
such that
\begin{equation}\label{nabla}
(\delta-\dbar)(v_1+v_2+\cdots)=\phi,
\end{equation}
i.e.,
\begin{equation}\label{nablaprim}
\delta v_{k+1}=\dbar v_k,\ k\ge 1,\quad \delta v_1(=fv_1)=\phi.
\end{equation}
For degree reasons, $\dbar v_k=0$ if $k$ is large enough, and if there are no obstructions
for solving $\dbar$, we can successively find $(0,k-2)$-forms (currents) $w_k$ with values
in $E_k$ such that 
\begin{equation}\label{losa}
\dbar w_k=v_k+\delta w_{k+1}, \quad    k\ge 2.
\end{equation}
Then finally
\begin{equation}\label{losning}
\psi=v_1+\delta w_2
\end{equation}
is a holomorphic solution to $f\psi=\phi$.
Since the $\dbar$-equations always are solvable locally we have

\begin{lma}\label{holo}
Suppose that we have a current solution $v=v_1+v_2+\ldots$ to  \eqref{nabla}. Then
 locally there are   holomorphic solutions  to $f\psi=\phi$.
\end{lma}

\smallskip


\smallskip
If $f$ is surjective, then obviously there are local holomorphic solutions to
$f\psi=\phi$ so the interesting case is when $f$ is just generically surjective.
In view of the argument in the proof, one also gets a global holomorphic  solution
provided all the $\dbar$-equations have global solutions.
Before we proceed with our construction  let us consider some examples.

\begin{ex}\label{stava}
If the complex \eqref{complex} is exact (in particular $f$ is surjective), 
then we can always
find such a solution to \eqref{nabla}.
In fact, given a holomorphic section $\phi$ of  $Q$, let $v_1$
be any pointwise solution to $\delta v=\phi$. Then
$\dbar v_1$ is $\delta$-exact and hence there is a $v_2$ such that 
$\delta v_2=\dbar v_1$ etc. 
\end{ex}

\begin{ex}
If $r=1$, i.e., $Q$ is a line bundle, then one can take $E_k=\Lambda^k E\otimes Q^*$  and 
$\delta$ as interior multiplication with $f$. One then gets  the usual
Koszul complex 
\begin{equation}\label{koszul}
\cdots\to \Lambda^3E\otimes (Q^*)^2\stackrel{\delta}{\to}\Lambda^2 E\otimes Q^*\stackrel{\delta}{\to}
 E\stackrel{f}{\to} Q\to 0,
\end{equation}
which is exact if (and only if) $f$ is non-vanishing.
 In fact, if we choose any section of  $\Hom(Q,E)\simeq E\otimes Q^*$
 $f\sigma=I_Q$ (if $E$ has a Hermitian metric  we can, e.g., choose
the section with pointwise minimal norm), then
there is  an induced mapping
$\sigma\colon \Lambda^{k}E\otimes (Q^*)^{k-1}\to \Lambda^{k+1}E\otimes (Q^*)^{k}$
such that and $\delta\circ \sigma+ \sigma\circ\delta=I$, 
and thus  \eqref{koszul} is exact. 
\end{ex}

\begin{ex}
Provided that $f$ is surjective, a simple way to find an exact complex \eqref{complex} is by taking
 $E_2 =\Ker f$ and  $E_k=0$ for $k>2$. 
However, $E_2$ is usually not trivial in  a neighborhood
of a singular point, i.e., $E_2$ usually cannot be extended as a vector bundle
across the set $Z$ where $f$ is not surjective.  To carry out the scheme in the proof of Lemma~\ref{holo}, one 
therefore 
has to solve a $\dbar$-equation in the bundle $E_2=\Ker f$ over  $X\setminus Z$. 
This is only possible under certain geometric conditions; this is precisely
what is  investigated and explained  in \cite{Sk2}.
\end{ex}

Let us now describe our  generalized Koszul complex.
Notice that our holomorphic morphism  $f\colon E\to Q$ is a holomorphic section
of the bundle  $\Hom(E,Q)$,  which we identify with
$E^*\otimes Q$.
If $\epsilon_j$ is a local holomorphic frame for $Q$, then
$$
f=\sum f_j\otimes\epsilon_j,
$$
where $f_j$ are sections of  $E^*$.
If $\eta$ is a section of $E$, then $f\eta=\sum_1^r\delta_{f_j}\epsilon_j$,
where $\delta_{f_j}$ denotes interior multiplication with $f_j$.
We can associate to $f$ the section
$$
F=f_1\w f_2\ldots\w f_r\otimes \epsilon_1\w\ldots\w\epsilon_r
$$
of $\Lambda^r E^*\otimes \det Q^*$.
It is   independent of the particular choice of frame, and  will be called  the
determinant section of $f$. Notice that $f$ is surjective at a point if and only of $F$ is non-vanishing
at that point. 
There is an induced mapping 
$$
\delta_F\colon \Lambda^{r+1}E\otimes\det Q^*\to E
$$
defined by
$$
\delta_F(\xi\otimes \epsilon_1^*\w\ldots\w\epsilon_r^*)=\delta_{f_r}\cdots\delta_{f_1}\xi,
$$
which is also  easily seen to be independent of the particular local frame $\epsilon_j$
for  $Q$; here $\epsilon_j^*$ denotes the dual frame for $Q^*$.
Moreover, it is also clear that
\begin{equation}\label{samm}
f\circ \delta_F=0.
\end{equation}

\smallskip

In order to proceed with the construction of our complex we have to recall some facts about
symmetric tensors.
Let $S^\ell Q^*$ be the subbundle of $\bigotimes Q^*$ consisting of 
symmetric $\ell$-tensors of $Q^*$. 
If $u,v\in Q^*$ then
$u\dot\otimes v=u\otimes v+v\otimes u$, etc. This 
extends to a commutative mapping $SQ^*\times SQ^*\to SQ^*$.
If $q$ is a section of  $Q$, then it induces the usual interior multiplication on
$\bigotimes Q^*$ (say from the left),  and,  in particular, if
$u^p=u\dot\otimes\cdots\dot\otimes u$, then
$\delta_q u^p=pu^{p-1} (q\cdot u)$.

\smallskip
We now define 
\begin{equation}\label{ekdef}
E_k=\Lambda^{r+k-1}E\otimes S^{k-2}Q^*\otimes\det Q^*, \quad k\ge 2.
\end{equation}
Given the local frame $\epsilon_j$, a section $\xi$ of  $E_k$ can be written
$$
\xi=\sum_{|\alpha|=k-2}\xi_\alpha\otimes \epsilon^*_\alpha\otimes \epsilon^*,
$$
where 
$$
\epsilon^*_\alpha=\frac{(\epsilon_1^*)^{\alpha_1}\dot\otimes\cdots\dot\otimes(\epsilon_r^*)^{\alpha_r}}
{\alpha_1!\cdots \alpha_r!}
$$
and
$
\epsilon^*=\epsilon_1^*\w\ldots\w\epsilon_r^*.
$
For $k\ge 2$ we  have  mappings
$
\delta\colon E_{k+1}\to E_k
$
defined by 
$$
\xi\otimes q^*\otimes \epsilon^*\mapsto
\sum_{j=1}^r \delta_{f_j}\xi\otimes \delta_{\epsilon_j} q^*\otimes \epsilon^*,
$$
which are also  independent of the specific choice of local frame $\epsilon_j$.
Since $\delta_{f_j}$ anti-commute and $\delta_{\epsilon_j}$ commute, it follows that 
$\delta^2=0$. Moreover, if the section $\xi$ of $E_2$ is in the image of $\delta$, i.e.,
 $\xi=\delta\eta=\sum \delta_{f_j}\eta_j\otimes\epsilon^*$, 
then clearly $\delta_F\xi=0$.  In view of \eqref{samm} we thus  have a complex
\begin{equation}\label{komplex}
\cdots\stackrel{\delta}{\to}E_3
\stackrel{\delta}{\to}E_2
\stackrel{\delta_F}{\to} E\stackrel{f}{\to} Q\to 0.
\end{equation}
In the sequel, we will often denote all the mappings in \eqref{komplex}
by $\delta$. 
Observe that if $r=1$, then \eqref{komplex} is just the Koszul complex
\eqref{koszul}.

\smallskip

If we let $\xi$  above take values in $\Lambda(T^*(X)_{0,1}\oplus E)$  rather than just $\Lambda E$, then
we get an extension of all the mappings $\delta$ and $\delta_F$ 
to forms and currents with values in $E_k$.
The mappings $\delta\colon E_{k+1}\to E_k$, $k\ge 2$, will automatically
satisfy \eqref{extend} so that $\delta\dbar=-\dbar\delta$, but we should 
have to  insert the factor $(-1)^{(r+1)q}$ in the definition
of $\delta_F$, when it acts on $\xi\otimes w$, and $\xi$ is a $(0,q)$-form.
However, we will not do that, and therefore we have instead that
$\dbar\delta_F=(-1)^r\delta_F\dbar$.
This means that the final solution $\psi$ in \eqref{losning} is
$\psi=v_1+(-1)^{r+1}\delta w_2$.

\smallskip

\section{Surjective morphisms} 

Now we  assume that $f\colon E\to Q$ is surjective and that $E$ and $Q$ are equipped with Hermitian
metrics. Moreover, we let $E_k$ be the the derived bundles defined by \eqref{ekdef}.
Let $\sigma$ be the section of  $\Hom(Q,E)=E\otimes Q^*$ with pointwise minimal norm 
(i.e., such that $\Im \sigma$ is orthogonal to $\Ker f$) such that 
$
f\circ\sigma=I_Q.
$
If $\epsilon_j^*$ denotes the dual frame for $Q^*$, then  
$$
\sigma=\sum\sigma_j\otimes\epsilon_j^*,
$$
where $\sigma_j$ are the sections of  $E$ with minimal  norms such that
$f_j\cdot \sigma_k=\delta_{jk}$.
Moreover, 
$$
\lsigma= \sigma_1\w\ldots\w\sigma_r\otimes \epsilon_1^*\w\ldots\w\epsilon_r^*
$$
is a well-defined section of $\Lambda E^r\otimes\det Q^*$, and it induces  a mapping
$E\to E_2=\Lambda^{r+1}E\otimes\det Q^*$ defined by
$$
\quad \xi\mapsto \lsigma\w\xi= \sigma_1\w\ldots\w\sigma_r\w \xi\otimes \epsilon_1^*\w\ldots\w\epsilon_r^*.
$$
Now, 
$\delta_F\lsigma\xi=\xi$ provided that $\delta_{f_j}\xi=0$ for all $j$, i.e.,
$\xi$ is in $\Ker f$.
Thus \eqref{komplex} is exact at $E$ if $f$ is surjective. We also have

\begin{lma}\label{snorkel}
If  $f$ is surjective, 
then \eqref{komplex} is  exact up to $E_2$.
\end{lma}

\begin{proof}
It remains to check the exactness at $E_2$.
Suppose that $\xi\otimes\epsilon^*$ is a section of  $E_2$ such that 
$0=\delta_F\xi\otimes\epsilon^*$.
Thus  $\xi$ is a section of 
 $\Lambda^{r+1}E$, such that $\delta_{f_r}\cdots\delta_{f_1}\xi=0$.
By the surjectivity of $f$, $\sigma_j$ are linearly independent and therefore they form a 
part of a basis $\sigma_1,\ldots,\sigma_m$ for $E$.
With respect to this basis we can write
$$
\xi=\sigma_1\w\ldots\w\sigma_r\wedge\xi'+\xi'',
$$
where $\xi''$ does not contain all the $\sigma_j$, $j=1,\ldots,r$. It follows that
$\xi'=0$ and hence
$\xi=\sum_1^r \xi_j$, where $\xi_j$ does  not contain $\sigma_j$.
Therefore
$$
\xi\otimes\epsilon^*=\sum_1^r\delta_{f_j}(\sigma_j\wedge\xi_j)\otimes\epsilon^*=
\delta\big(\sum_1^r\sigma_j\w\xi_j\otimes\epsilon^*_j\otimes\epsilon^*\big),
$$
and thus \eqref{komplex} is exact at $E_2$.
\end{proof}


In order to find a solution to \eqref{nabla} it is natural to start with
$v_1=\sigma\phi=\sum \phi_j\sigma_j$, and
$v_2=\lsigma\w \dbar\sigma\phi$.
We can just as well suppress $\phi$ and define $u_k$ with values
in 
$\Hom(Q,E_k)$ such that $u_k\phi$ satisfies \eqref{nabla}. 
Notice that $\Hom(Q,E_1)\simeq E\otimes Q^*$ and
$$
\Hom(Q,E_k)\simeq \Lambda^{r+k-1}E\otimes S^{k-2}Q^*\otimes\det Q^*\otimes Q^*, \quad k\ge 2.
$$
\begin{df}
We define the $(0,k-1)$-forms $u_k$ with values in $\Hom(Q,E_k)$ as
\begin{equation}\label{udef1}
u_1=\sigma, \quad
u_k=(\dbar\sigma)^{\otimes(k-2)}\otimes\lsigma\otimes\dbar\sigma, \quad k\ge 2.
\end{equation}
\end{df}

Notice that $(\dbar\sigma)^{\otimes(k-2)}$ is indeed a symmetric tensor,
and that the definition of $u_k$ is invariant.
Moreover, 
$$
{\big(\sum_1^r\dbar\sigma_\ell\otimes\epsilon^*_\ell\big)^{\otimes(k-2)}}=
\frac{(\sum_1^r\dbar\sigma_\ell\otimes\epsilon^*_\ell)^{\dot\otimes(k-2)}}
{(k-2)!}=
\sum_{|\alpha|=k-2}(\dbar\sigma)^\alpha
\otimes
\epsilon_\alpha^*,
$$
where 
$$
(\dbar\sigma)^\alpha=(\dbar\sigma_1)^{\alpha_1}\w\ldots\w(\dbar\sigma_r)^{\alpha_r}.
$$
Since
$\dbar\sigma_j$ have degree $2$ in $\Lambda(T^*_{0,1}(X)\oplus E)$ and therefore
commute, we thus have that 
\begin{multline}\label{udef}
u_1=\sigma=\sum_j\sigma_j\otimes\epsilon_j^*, \\
u_k=\sigma_1\w\ldots\w\sigma_r\w\sum_{|\alpha|=k-2}\sum_j(\dbar\sigma)^\alpha
\w\dbar\sigma_j\otimes
\epsilon_\alpha^*\otimes\epsilon^*\otimes\epsilon_j^*,\quad k\ge 2.
\end{multline}

\begin{prop}
Assume that $f$ is surjective.
If $u$ is defined by \eqref{udef1}, then
$$
(\delta-\dbar)(u_1+u_2+\cdots )=I_Q.
$$
\end{prop}

Here $I_Q\colon Q\to Q$ is the the identity morphism.

\begin{proof}
We have already seen that $\delta u_1= f\sigma=I_Q$ and $\delta_F\lsigma\otimes\dbar\sigma=
\dbar\sigma=\dbar u_1$, so
we  have to verify that $\delta u_{k+1}=\dbar u_k$ for $k\ge 2$.  Now,
$$
u_{k+1}=
\sigma_1\w\ldots\w\sigma_r\w {\big(\sum_1^r\dbar\sigma_\ell\otimes\epsilon^*_\ell\big)^{\otimes(k-1)}}
\otimes\epsilon^*\otimes\dbar\sigma.
$$
Recalling  that  $\delta=\sum\delta_{f_j}\otimes\delta_{\epsilon_j}$, and that 
$\delta_{\epsilon_j}$ acts from the  left,  we get
$$
\delta u_{k+1}=\sum_{j=1}^r\delta_{f_j}(\sigma_1\w\ldots\w\sigma_r)\w \dbar\sigma_j\w
(\sum_1^r\dbar\sigma_\ell\otimes\epsilon^*_\ell)^{\otimes(k-2)}\otimes\epsilon^*\otimes\dbar\sigma.
$$
On the other hand,
$$
\dbar u_k=
\dbar(\sigma_1\w\ldots\w\sigma_r)\w (\sum_1^r\dbar\sigma_\ell\otimes\epsilon^*_\ell)^{\otimes(k-2)}
\otimes\epsilon^*\otimes\dbar\sigma.
$$
It is now clear that $\delta u_{k+1}=\dbar u_k$.
\end{proof}

It follows that if  $\phi=\sum\phi_j\epsilon_j$ is a holomorphic section of  $Q$ and $v_k=u_k\phi$, then
$v=u\phi$ satisfies \eqref{nabla}.
For later purpose  we rewrite the expression for $u\phi$ so that it only
``depends'' on  $\sigma\phi=\sum_j \phi_j\sigma_j$. In fact, applying
$\dbar$ to the equality 
$0=\lsigma\otimes\sigma$ we get
$0=\lsigma\otimes\dbar\sigma+(-1)^r\dbar\lsigma\otimes\sigma$.
Therefore,
\begin{equation}\label{omskriv}
u_k\phi=(-1)^{r+1}(\dbar\sigma)^{\otimes(k-2)}\otimes\dbar\lsigma\otimes\sigma\phi, \quad k\ge 2,
\end{equation}
or more explicitly (the possible minus sign cancels out)
\begin{multline*}
 u_k\phi=\\
(\sum_j\phi_j\sigma_j)\w\dbar(\sigma_1\w\ldots\w\sigma_r)\w
\sum_{|\alpha|=k-2}
(\dbar\sigma)^\alpha
\otimes
\epsilon_\alpha^*\otimes\epsilon^*=\\
(\sum_j\phi_j\sigma_j)\w\dbar(\sigma_1\w\ldots\w\sigma_r)
\w(\sum_1^r\dbar\sigma_\ell\otimes\epsilon^*_\ell)^{\otimes(k-2)}
\otimes\epsilon^*,
\quad k\ge 2.
\end{multline*}

\smallskip

\section{Analysis of  singularities}

We now consider an $f$ that is not necessarily surjective everywhere. Then
we can define $u$ as above in $X\setminus Z$, where
$Z$ is the set where $f$ is not surjective, which
is equal to the zero set of the holomorphic section $F$
and hence an analytic subvariety of $X$.
To analyze the singularities of $u$ at $Z$ we will use the following lemma.

\begin{lma}\label{fakt}
(i)\  Let $s$ be the section of  $\Hom(Q,E)\simeq E\otimes Q^*$ 
with pointwise minimal norm such that 
$$
fs=|F|^2 I_Q,
$$
and let $S$ be the section of $\Lambda^rE\otimes\det Q^*$ with pointwise minimal norm such that
$$
FS=1.
$$
Then
$s$ and $S$ are smooth across $Z$, and 
$$
s= |F|^2\sigma,\quad  S=|F|^2\lsigma   \quad {\text in} \ X\setminus Z.
$$
\smallskip

\noindent (ii) \  If in addition  $F=F_0 F'$ for some holomorphic  function $F_0$ and 
non-vanishing holomorphic section $F'$, then 
$$
s'=F_0\sigma, \quad S'=F_0\lsigma
$$
are smooth across  $Z$.
\end{lma}

Given a section $\eta$ of $E^*$, let $\eta^*$ denote the dual section with respect to 
the Hermitian metric, i.e., $\langle \xi, \eta^*\rangle =\eta\cdot\xi$
for sections $\xi$ of $E$.  The mapping $\eta\mapsto \eta^*$ is conjugate-linear
and extends to a conjugate-linear mapping
$\Lambda E^*\to\Lambda E$ by 
$$
\eta_1\w\ldots\w\eta_r\mapsto \eta_1^*\w\ldots\w\eta_r^*.
$$

\begin{proof}
Assume that $\epsilon_j$ is a local frame for $Q$ as before,  let $\epsilon_j^*$ be its
dual frame, and assume $f=\sum f_j\otimes\epsilon_j$. Now, 
$$
S=|\epsilon_1\w\ldots\w\epsilon_r|^2 f_1^*\w\ldots\w f_r^*\otimes\epsilon_1^*\w\ldots\w\epsilon_r^*
$$
is a section of $\Lambda^r E\otimes\det Q^*$ that is independent of the particular choice
of frame. This is checked by considering a change of frame
$\epsilon'=\epsilon g$, where $g$ is an invertible $r\times r$-matrix. 
Notice that
\begin{equation}\label{f2}
FS=|\epsilon_1\w\ldots\w\epsilon_r|^2|f_1\w\ldots\w f_r|^2=|F|^2.
\end{equation}
We can choose the frame $\epsilon_j$ such that   $f_j$ are orthogonal at any given point
outside $Z$. Thus  
 $f_j=\alpha_je_j^*$, $j=1,\ldots,r$ for some $ON$-frame  $e_j^*$ of $E^*$ and 
$f_j^*=\bar\alpha_j e_j$, and it  is then easy to see that $S$ is in fact the section with minimal
norm such that $FS=1$. 
In particular this means that $S$ is the dual section of $F$.
Moreover, at this point, $\sigma_j=(1/\alpha_j)e_j$, $j=1,\ldots,r$, and thus
$$
\lsigma=\frac{e_1\w\ldots\w e_r}{\alpha_1\cdots\alpha_r}\otimes\epsilon_1^*\w\ldots\w\epsilon_r^*.
$$
Therefore,
$$
|F|^2\lsigma=|\epsilon_1\w\ldots\w\epsilon_r|^2\bar\alpha_1\cdots\bar\alpha_r e_1\w\ldots\w e_r
\otimes \epsilon_1^*\w\ldots\w\epsilon_r^*
=S.
$$
Now assume that  $F=F_0F'$. Since $S$ is the dual section of $F$ it follows that
$S=\bar F_0S''$, where $S''$ is the dual section of $F'$, and thus $S''$ is smooth even
across $Z$.
Therefore  $F_0\lsigma=F_0 S/|F|^2=S''/|F'|^2$ is smooth across $Z$ as well.

Now, define
\begin{multline*}
s=\big(\sum_1^r\delta_{f_j}\otimes\delta_{\epsilon_j}\big)^{r-1}S/r!= \\
|\epsilon_1\w\ldots\w\epsilon_r|^2
\sum_1^r(-1)^{\ell+1}\delta_{f_r}\cdots\delta_{f_{\ell+1}}\delta_{f_{\ell-1}}\cdots\delta_{f_1}
(f_1^*\w\ldots\w f_r^*)\otimes\epsilon_j^*.
\end{multline*}
Clearly,
$$
fs=\sum_1^r\delta_{f_j}\otimes{\epsilon_j}s=|\epsilon_1\w\ldots\w\epsilon_r|^2|f_1\w\ldots\w f_r|^2\sum_j\epsilon_j
\otimes\epsilon_j^*=|F|^2 I_Q.
$$
Moreover, it is readily checked that $s$ is orthogonal to $\Ker f$ so that $s$ is the minimal section
such that $fs=|F|^2 I_Q$. One can also check this by choosing a frame as above such that  that $f_j$ are
orthogonal. Thus, $s/|f|^2=\sigma$.
Finally, if $F_0\lsigma$ is smooth across $Z$ it follows that
$$
F_0\sigma=\big(\sum_1^r\delta_{f_j}\otimes\delta_{\epsilon_j}\big)^{r-1}F_0\lsigma
$$
is smooth as well.
\end{proof}

\begin{remark}\label{travat}
Assume that $Q$ as well as $E$ are trivial and equipped with the trivial metrics,
and  $e_j$ and $\epsilon_j$  are ON-frames. If 
$$
F=\sum'_{|I|=r}F_I \otimes e^*_I,
$$
(suppressing the factor $\epsilon_1\w\ldots\w\epsilon_r$ and its dual), 
we have that
$$
S=\sum'_{|I|=r}\overline{F_I} \otimes e_I.
$$
\end{remark}

Since $\lsigma\otimes\sigma=0$ we have, for any scalar function $\xi$,  that
\begin{multline}\label{homo}
(\dbar(\xi\sigma))^{\otimes(k-2)}\otimes \xi \lsigma\otimes\dbar(\xi\sigma)=\\
(-1)^{r+1}
(\dbar(\xi\sigma))^{\otimes(k-2)}\otimes \dbar(\xi \lsigma)\otimes \xi\sigma=
\xi^k u_k.
\end{multline}
In view of  Lemma~\ref{fakt} we have in  particular  that 
\begin{multline}\label{nyform}
u_1=\frac{s}{|F|^2}, \\
u_k=\frac{(\dbar s)^{\otimes(k-2)}\otimes\ S\otimes \dbar s}{|F|^{2k}}
=(-1)^{r+1}\frac{(\dbar s)^{\otimes(k-2)}\otimes\ \dbar S\otimes s}{|F|^{2k}}, \quad k\ge 2
\end{multline}


\begin{lma} \label{normus}
If  $f^*\colon Q\to E$ is the adjoint morphism with respect to  the Hermitian structures
on $Q$ and $E$,  then 
\begin{equation}\label{bulle}
|F|^2=\det (ff^*).
\end{equation}
\end{lma}

\begin{proof}
Since both sides of \eqref{bulle} are invariant pointwise statements, we can assume that $\epsilon_j$ is a
$ON$-frame. Let $\xi$ be any section of $E$. Then
$$ 
f_\ell\cdot\xi=\langle\sum_j(f_j\cdot\xi)\epsilon_j,\epsilon_\ell\rangle =
\langle f\xi,\epsilon_\ell\rangle =\langle\xi, f^*\epsilon_\ell\rangle
$$
which means that $f^*=\sum_j f_j^*\otimes \epsilon_j^*$.
It follows that
$$
ff^*=\sum (f_j\cdot f_k^*) \epsilon_j\otimes\epsilon^*_k=\langle f_j,f_k\rangle \epsilon_j\otimes\epsilon^*_k.
$$
Thus 
$$
\det (ff^*)=\det \langle f_j,f_k\rangle = |f_1\w\ldots\w f_r|^2,
$$
which implies the statement since $|\epsilon_1\w\ldots\w \epsilon_r|=1$, cf., \eqref{f2}.
\end{proof}

Since  $(ff^*)^{-1}=\widetilde{ff^*}/\det(ff^*)$ and
$s/|F|^2=\sigma=f^*(ff^*)^{-1}$ we can conclude that
$$
s=f^*\widetilde{ff^*}.
$$
Moreover, if $F=F_0F'$ as in Lemma~\ref{fakt} above, then by \eqref{bulle}, 
\begin{equation}\label{normis}
|s'\phi|=|\sigma\phi||F_0|=|f^*(ff^*)^{-1}\phi||F|/|F'|=
\|\phi\|/|F'|.
\end{equation}

\vspace{.5cm}

\section{The residue current of a generically surjective holomorphic morphism}\label{rescurrent}

We say that $f\colon E\to Q$ is generically surjective if $Z$ has positive codimension. 
Again let $E_k$ be the derived bundles defined by \eqref{ekdef} and let
$H=Q\oplus E\oplus E_2\oplus E_3\oplus\cdots$.
From Section~3  we have we the section  $u$  of 
$\Hom(Q,H)$ over $X\setminus Z$   defined by \eqref{udef}.
Following \cite{A2} we shall now extend it to a current
$U$ with values in $\Hom(Q,H)$ across $Z$ and  define the corresponding residue current.

\begin{thm}\label{bas}
Assume that $f\colon E\to Q$ is a generically surjective holomorphic morphism and let
$u$ be the associated section of $\Hom(Q,H)$ defined in $X\setminus Z$. The function
$\lambda\mapsto |F|^{2\lambda}u$ is holomorphic for $\Re\lambda>-\epsilon$ and
$$
U= |F|^{2\lambda}u|_{\lambda=0}
$$
 is a current  extension of $u$ across $Z$. Moreover, 
$$
(\delta-\dbar) U=I_Q-R^f,
$$
where
$$
R^f= \dbar |F|^{2\lambda}\w u|_{\lambda=0}.
$$
The current $R^f$ (taking values  in $\Hom(Q,H)$) has support on $Z$ and 
$$
R^f=R^f_p+\cdots +R^f_\mu,
$$
where $p=\codim Z$ and $\mu=\min(n,m-r+1)$,
and $R^f_k$ is a current of bidegree $(0,k)$ with values in $E_k$.  
\end{thm}

\begin{proof}
The proof is more or less identical to the proof of Theorem~1.1 in \cite{A2}, so we only
sketch it. 
After  an appropriate  resolution of singularities we may assume  that 
$F=F_0 F'$, where $F_0$ is a holomorphic  function and $F'$ is a non-vanishing
section of  (the pullback of) $\Lambda^r E^*\otimes\det Q$. 
According to Lemma~\ref{fakt}~(ii), then  $s'=F_0\sigma$ and $S'=F_0\lsigma$ are  smooth across the singularity, and
hence by \eqref{homo}, 
$$
u_1=\frac{s'}{F_0},\quad
u_k=(-1)^{r+1} \frac{(\dbar s')^{\otimes(k-2)}\otimes\dbar S'\otimes s'}{F_0^k},\quad k\ge 2.
$$
It is now easy to see that the analytic continuations
of $|F|^{2\lambda}u_k$ exist, and 
in this resolution the values at $\lambda=0$ are  just the  principal value currents
$$
u_1=\Big[\frac{1}{F_0}\Big]s',\quad  \Big[\frac{1}{F_0^k}\Big](\dbar s')^{\otimes(k-2)}\otimes\dbar S'\otimes s'.
$$
Precisely as in \cite{A2}, by the way following \cite{BY1} and  \cite{PTY}, one can show  that 
$$
R_k^f=\dbar|F|^{2\lambda}\w u_k\big|_{\lambda=0}=0
$$ 
if $k<p=\codim Z$. Thus Theorem~\ref{bas} is proved.
\end{proof}

\begin{proof}[Proof of Theorem~\ref{omv}]
If 
$$
\phi= f\psi=\sum_1^r (\delta_{f_j}\psi) \epsilon_j, 
$$
then 
\begin{multline*}
u_{m-r+1}\phi=(\sum_1^r\phi_j\dbar\sigma_j)\w
\sigma_1\w\ldots\w\sigma_r\w(\dbar\sigma)^{\otimes(m-r-1)}
\otimes\epsilon^* =\\
(\sum_{j=1}^r\delta_{f_j}\psi\w\dbar\sigma_j)\w
\sigma_1\w\ldots\w\sigma_r\w(\dbar\sigma)^{\otimes(m-r-1)}
\otimes\epsilon^*=\\
\psi\w\sum_{j=1}^r\delta_{f_j}\big(\sigma_1\w\ldots\w\sigma_r\big)\w\dbar\sigma_j
\w(\dbar\sigma)^{\otimes(m-r-1)}\otimes
\epsilon^*=\\
\dbar\Big(\psi\w\sigma_1\w\ldots\w\sigma_r
\w(\dbar\sigma)^{\otimes(m-r-1)}
\otimes\epsilon^*\Big) =\dbar (u'\psi),
\end{multline*}
where we have used that the form has maximal  degree $m$  in $\Lambda E$.
If we define
$$
R'\psi=\dbar|f|^{2\lambda}\w u'\w\psi|_{\lambda=0},
$$
it follows that 
$$
R_{m-r+1}\phi=\dbar (R'\psi).
$$
However, since  $\codim Z=m-r+1$ it follows as in the proof of Theorem~\ref{bas} above that
$R'\psi$,  being a $(0,m-r)$-current, vanishes for degree reasons.
Thus the theorem is proved.
\end{proof}


\begin{proof}[Proof of Proposition~\ref{bsgen}]
Let $\xi$ by a test form with support contained in neighborhood where we have the
resolution of singularities.
In view of \eqref{homo} we see that
$R^f_k\cdot\xi$ is a sum of terms like
$$
\dbar|F_0|^{2\lambda}v^\lambda\w
\frac{(\dbar s')^{\otimes(k-2)}\otimes\dbar S'\otimes s'\phi}
{F_0^{k}}\rho\Big|_{\lambda=0},
$$
where
$\rho$ is a cut-off function, and $v$ is a smooth strictly positive function.
By the hypothesis, \eqref{normis}, and \eqref{normus}, 
$$
|s'\phi|\sim \|\phi\|\lesssim |F|^{\min(n,m-r+1)}\lesssim
|F_0|^{k},
$$
since $k\le\min(n,m-r+1)$. Thus we  must check that
$$
\int\dbar|F_0|^{2\lambda}v^\lambda\w\frac{\eta}{F_0^\ell}
$$
vanishes at $\lambda=0$ if $v$ is a test form that is 
$\O(|F_0|^\ell)$. However, we may as well assume that $F_0$ is a monomial in the
local coordinates, and therefore this statement follows from the corresponding one-variable
statement;  that
$$
\int_\tau\dbar|\tau|^{2\lambda}v^\lambda\w \frac{\eta(\tau)d\tau}{\tau^\ell}
$$
vanishes at $\lambda=0$ if $\eta=\O(|\tau|^\ell)$.  If we write
$$
\dbar|\tau|^{2\lambda}=\lambda|\tau|^{2\lambda}\frac{d\bar\tau}{\bar\tau}
$$
this follows by dominated convergence.
\end{proof}

\begin{proof}[Proofs of Theorems~\ref{basal2} and \ref{bsgen2}]
These theorems  are proved in much the same way as the corresponding results
for holomorphic functions.
In fact, if $\phi$ is a $\dbar$-closed smooth form with values in $Q$, then
$(\delta-\dbar) U\phi =\phi-R^f\phi$
so  $(\delta-\dbar)U\phi=\phi$ if $R^f\phi=0$. Following an apparent modification
of the   procedure
in \eqref{losa}  we get a $\dbar$-closed current $\psi$ with values in $E$ such that
$f\psi=\phi$. 
However, if  \eqref{bsvillkor2} holds, then it is not hard to
verify that $U\phi$ is locally integrable; in fact in a resolution
as above $U\phi$ is then bounded, in particular it is locally integrable,
and therefore $U\phi$ is locally integrable in $X$.
 Therefore, we can get a locally integrable solution  $\psi$.
\end{proof}

For further reference we state the following  proposition  that shows  that the 
principal term  $R_p^f$ of the current  $R^f$ (where $p=\codim Z)$ is  robust
in certain sense. It is proved precisely as Theorem~2.1 in \cite{A3}.

\begin{prop}\label{gulnad}
Let $f$ be a generically surjective morphism, let $R^f$ be the associated
residue current, and let $p=\codim Z$. 
Assume that $h$ is a holomorphic section of some line bundle such that  
$\{h=0\}\cap Z$ has codimension $p+1$. If $\phi$ is a holomorphic section such that
 $R_p^f\phi=0$ in $X\setminus\{h=0\}$, then 
$R^f_p\phi=0$ in $X$.
\end{prop}

\vspace{.5cm}

\section{Polynomial mappings}

Let $L^s=\O(s)$ be the line bundle over $\P^n$ whose sections are represented by
$s$-homogenoues functions in $z=(z_0,\ldots,z_n)$.
If $\phi$ is a section of $L^s$ its natural norm is
$$
\|\phi(z)\|=\frac{|\phi(z)|}{|z|^s}.
$$
Let $E_j$, $j=1,\ldots,m$,  be trivial line bundles over $\P^n$ with basis
elements $e_j$. If $f^j$ are the $d_j$-homogenizations of the columns of polynomials
$P^j$ as in Theorem~\ref{pellenom}, then
$f_k=\sum_1^m f_k^j e_j^*$ are sections of 
$$
E^*=E_1^*\otimes L^{d_1}\oplus\cdots\oplus E^*_m\otimes L^{d_m}.
$$
Observe that  the section $F=f_1\ldots\w f_r$  of $\Lambda^r E^*$ 
can be written
$$
F=\sum_{|I|=r}' F_I e^*_{I_1}\w\ldots\w e^*_{I_r},
$$
where
$
F_I=\det(f_k^{I_j}).
$
Thus
$$
\|F(z)\|^2=\frac{\sum_{|I|=r}'|F_I(z)|^2}{|z|^{2\sum_1^rd_{I_j}}}.
$$
If we write this expression in the affine coordinates $z'$ we get precisely the left hand
side of \eqref{olik}.
Let  $Q$ be the trivial bundle $\C^r\to\P^n$. Then $f$ defines a morphism
 $f\colon E\to Q$ such that 
$$
\psi=(\psi_1,\ldots,\psi_m)\mapsto
\sum_1^m f^j\psi_j
$$
 where $f^j$ are the homogenizations of the given
(columns) of polynomials $P_j$ of degrees $d_j$.
We can now prove  Theorem~\ref{pellenom}.

\begin{proof}[Proof  of Theorem~\ref{pellenom}]
Let
$$
\cdots \to E_2\to E\to \C^r\to 0.
$$
be the induced complex defined in Section~\ref{nykomplex}. We can take 
tensor products with $L^\rho$ and get  the complex
$$
\cdots \to E_2\otimes L^\rho\to E\otimes L^\rho\to \C^r\otimes L^\rho\to 0.
$$
Let $U$ and $R^f$ be the corresponding currents from
Section~\ref{rescurrent} with respect to the natural Hermitian metrics of $E$ and $Q$. If $\phi$ is a section of
$Q\otimes L^\rho$, and $R^f\phi=0$, then
$v=U\phi$ solves the equations
$$
(\delta-\dbar)(v_1+v_2+\cdots +v_{\min(n+1,m-r+1)})=\phi.
$$
In order to get a holomorphic solution $\psi$ to $f\psi=\phi$ we have to solve
all the equations $\dbar w_k= v_k-\delta v_{k+1}$.
Notice that  $v_k-\delta v_{k+1}$ is a $(0,k-1)$-current with values in  $E_k\otimes L^\rho$.
It is well-known that 
$H^{0,k}(\P^n,L^\nu)=0$ for  all $\nu$ if $1\le k\le n-1$, whereas 
$H^{0,n}(\P^n,L^\nu)=0$ (if and only) if $\nu\ge -n$.
If $m-r+1\le n$ there is therefore no problems at all, and 
the only possible obstruction may  appear  when $k=n+1$. 
Notice that  $v_{n+1}$ takes values in
$$
E_{n+1}\otimes L^\rho=\Lambda^{r+n}E\otimes S^{n-1}Q^*\otimes Q^*\otimes L^\rho.
$$
Since $Q$ is trivial, $E_{n+1}$ is a direct sum of line bundles 
$$
L^{\rho-(d_{I_1}+\ldots +d_{I_{r+n}})},
$$
where $I$ is an increasing multiindex. Therefore the crucial $\dbar$-equation
is solvable if $\rho-(d_1+\cdots +d_{n+r})\ge -n$.
Finally we express the relation 
$\sum f^j\psi_j=f\psi=\phi$ in affine coordinates and get the desired polynomials
$Q_j$ as $Q_j(z')=\psi_j(1,z')$.
\end{proof}

\begin{proof}[Proof of Proposition~\ref{kollest}]
If we have $m'$ polynomials $P'_\nu$ of degrees at most $d'$ with no common zeros in $\C^n$,
then it is proved in \cite{Koll} that 
$$
\frac {|P_\nu(z')|^2}{(1+|z'|^2)^{d'}}\ge C\frac{1}{(1+|z'|^2)^{M}},
$$
where $M=(d')^{\min(n,m')}$, provided that $d'\ge 3$.
We have $m!/(m-r)!r!$ polynomials $P_I=\det(P^k_{I_j})$ of degrees $d'=rd$ and hence
the proposition follows.
\end{proof}

\begin{proof}[Proof of  Theorem~\ref{noethergen}]
Since $\phi=f\psi$ is solvable in $\C^n$ and $p=m-r+1$,
$R^f\phi=0$ in $\C^n$ by Theorem~\ref{omv}.
If we take the section  $h=z_0$ of $\O(1)$, then
since assumption $Z$ has no irreducible component in 
the plane at infinity, 
$\codim Z\cap\{z_0=0\}=m-r+2$. By Proposition~\ref{gulnad}
we therefore have that $R^f\phi=0$ in $\P^n$.
Since \eqref{villkor} is fulfilled, the desired solution is given by
Theorem~\ref{pellenom}.
\end{proof}

\section{Estimates for a pointwise surjective morphism}\label{hpcorona}

In this section we indicate that our method can be used to get new quite sharp
estimates even when $f$ is pointwise surjective.
Let us assume that $E\simeq\C^m$ and $Q\simeq\C^r$ are trivial bundles over a
smoothly bounded domain $D=\{\rho<0\}$ in $\C^n$, and equipped with the trivial
metrics. Then a morphism $f\colon E\to Q$ is just a matrix of holomorphic functions in $D$.
We assume that $f\in H^\infty(D,\Hom(E,Q))$, and that moreover
$$
|f_1\w\ldots \w f_r|\ge \delta>0;
$$
this means that $f$ is uniformly surjective. 
Notice that since $Q$ is trivial, $\det Q^*$ is just the trivial line bundle, so
$F=f_1\w\ldots \w f_r$.
Let
$$
\|\phi\|^p_{H^p}=\sup_{\epsilon>0}\int_{\rho=-\epsilon}|\phi(z)|^pdS, \quad 0<p<\infty.
$$
It was proved in \cite{A7} and  \cite{A8}  that if $D$ admits a plursubharmonic  defining function
$\rho$,  and $p\le 2$, then for any $\phi\in H^p(D,Q)$
there is a $\psi\in H^p(E)$ such that $f\psi=\phi$. Morerover, the norm of $\phi$ is
bounded by a constant times $\log(1/\delta)/\delta^{1+\mu}$, where $\mu=\min(n,m-r)$.
To be precise, this sharp estimate is only explicitly given in the case $r=1$, but it 
follows (hopefully) in the general case as well with a similar argument. 
This result is proved by a combination of the $L^2$-methods in \cite{Sk2}, the refined $L^2$ estimate
for $\dbar_b$ 
introduced in \cite{B4},  and  Wolff type estimates.  
\smallskip

The case $p>2$ and $r=1$ has been studied by several authors in strictly pseudoconvex domains,
e.g., \cite{Am} and \cite{AC},
and a generalization to $r>1$ is made  in \cite{JH1}. These works are based on
integral representation and harmonic analysis. There are also similar results
for other spaces of functions, see, e.g., the references in \cite{AC}.
To show how the ideas in this paper can be applied in a situation like this we present
the new result Theorem~\ref{bollar} below.  It is clear that other known results that are proved
by means of the Koszul complex in the case $r=1$ can be 
generalized to the case $r>1$ in an analogous way.

\begin{thm}\label{bollar}
Let $D$ be a strictly pseudoconvex  domain with reasonably smooth  ($C^3$ is enough) 
boundary and $p<\infty$. For any $\phi\in  H^p(D,Q)$ there is a $\psi\in  H^p(D,E)$
such that $f\psi=\phi$ and
\begin{equation}\label{gaphals}
\|\psi\|_{H^p}\le C_\delta\|\phi\|_{H^p},
\end{equation}
where
$$
C_\delta\le C(\log(1/\delta))^{\mu/2}/\delta^{1+\mu},
\quad {\rm if}\quad \mu=\min(n,m-r)>1
$$
and
$$
C_\delta\le C\log(1/\delta)/\delta^{2}\quad {\rm if}\quad
\min(n,m-r)=1.
$$
and $C_\delta=1/\delta$ if $\min(n,m-r)=0$, i.e., $m=r$,
and $C$ is a constant that is independent of $m$. 
\end{thm}

In the case $r=1$ this coincides with the result in \cite{AC}.
Since the proof of this generalization follows the proof in
\cite{AC} quite closely we just give a sketch and indicate
the necessary modifications.

\begin{proof}[Sketch of proof]
With the notation as before we have that $v=u\phi$, where
$u=u_1+\cdots u_{1+\mu}$
and 
$$
u_{k+1}=(\dbar\sigma)^{\otimes(k-1)}\otimes\lsigma\otimes\dbar\sigma.
$$
Since $Q$ is trivial and $\epsilon_j$ is a ON-frame 
we have, cf.,  Remark~\ref{travat}, 
that $S=\sum' \overline{F_I} e_I\otimes \epsilon^*$ of
$F=\sum' F_I e^*_I\otimes\epsilon$, so the coefficients in $S$ are anti-holomorphic. 
 Moreover, $\lsigma = S/|F|^2$ and 
$$
\sigma=(\sum\delta_{f_j}\otimes\delta_{\epsilon_j})\lsigma/k!
$$
so that 
$$
\dbar\sigma=\pm(\sum\delta_{f_j}\otimes\delta_{\epsilon_j})\dbar\lsigma/k!.
$$
The norm of forms will be taken with respect to the non-isotropic
metric 
$\Omega=-\rho i\partial\dbar\log(1/-\rho))$; we assume that
$\rho$ is a strictly plurisubharmonic defining function.
For  $k\ge 2$ we can estimate the Carleson norm of 
$
(-\rho)^{\frac{k-2}{2}}|u_{k+1}|
$ 
precisely as in the proof of  Proposition~5.2 in \cite{AC}, with $F$ instead of $g$ and
$\lsigma$ instead of $\gamma$, and obtain the same estimate
$C\delta^{-k-1}(\log(1/\delta))^{k/2}$.

The necessary estimate of $u_2$ is of Wolff type and  involves
$\L u_2$ where $\L$ is a smooth $(1,0)$-vector field. More precisely
we  need to know that the Carleson norm of $|u_2|^2$
is bounded by $C^2$ and the Carleson norm of
$\sqrt{-\rho}|\L u_2|$ is bounded by $C$, where
$C=\delta^{-2}\log(1/\delta)$.
It is easily seen that 
$$
|u_2|\lesssim\frac{|\partial F|}{|F|^3}\le\frac{1}{\delta}\frac{|\partial F|}{|F|^2},
$$
and the desired Carleson estimate of $|u_2|^2$ now follows form Proposition~5.2 in \cite{AC}.
For any holomorphic function $\psi$ we have that $\sqrt{-\rho}|\L\psi|\lesssim||partial\psi|$.
When we compute $\L u_2$ we get either derivatives on the factor $1/|F|^2$ or on the
functions $f_j$ in $(\sum\delta_{f_j}\otimes\delta_{\epsilon_j})^{r-1}$. 
Therefore we have that
$$
\sqrt{-\rho}|\L u_2|\lesssim\frac{\partial h||\partial F|}{|F|^3}+\frac{|\partial F|^2}{|F|^4}
\lesssim \frac{1}{\delta^2}+\frac{|\partial F|^2}{|F|^4},
$$
where $h$ is a  holomorphic and bounded. It now follows from 
 Proposition~5.2 in \cite{AC} that the Carleson norm is $\lesssim \delta^{-2}\log(1/\delta)$
as wanted. 
\end{proof}

\begin{remark}
If we just assume that $\phi$ is in $H^p$ with respect to the stronger pointwise norm
$\|\phi\|$ instead of
 $|\phi|$,  the same proof would give a solution in $H^p$ provided that we could
say that the maximal function of $\|\phi\|$ is in $L^p(\partial D)$.
However, we do not know if this is true.
\end{remark}

As mentioned above, a similar result   has  previously been obtained by Hergoualch, \cite{JH1},
using an  idea of Fuhrmann, \cite{Fu}, to reduce to the case  $r=1$. 
However, by this method one has some loss of precision in the dependence of $\delta$.  
To see this let us  describe this method. Assume that $\phi=\sum\phi_j\epsilon_j$
is in $H^p(D,Q)$. 
Now  $F=f_1\w\ldots\w f_r\in  H^\infty(D,\Lambda^rE^*)$ and
$|F|\ge\delta$, so by the corresponding result for $r=1$, for each $\phi_j$ we can find
a section $H^j$ of $H^p(D,\Lambda^r E)$ such that $F\cdot H^j=\phi_j$, $j=1,\ldots,r$, with
$\|H^j\|_{H^p}\le C_{\delta}'\|\phi\|_{H^p}$.
 Since the rank of $\Lambda^r E^*$ is $m!/(m-r)!r!$, we get
$$
C'_\delta\le C\frac{(\log(1/\delta))^{\mu' /2}}{\delta^{1+\mu'}},
\quad \mu'=\min(n,m!/(m-r)!r!),
$$
if $\mu'>1$.
Since $f_1\w\ldots\w f_r\cdot H^j=\phi_j$, we also have that
$
f_j\cdot \psi_j=\phi_j
$
if 
$$
\psi_j= (-1)^{j+1} \delta_{f_r}\cdots \delta_{f_{j+1}}\delta_{f_{j-1}}\cdots\delta_{f_1} H^j.
$$
Now  $\psi=\psi_1+\cdots +\psi_r$ solves $f_j\psi=\phi_j$ for each $j$, i.e.,
$f\psi=\phi$ as wanted, and since $f_j$ are bounded,
we get the same estimate 
$$
\|\psi\|_{H^p}\le C'_{\delta}\|\phi\|_{H^p}.
$$
However, $\mu' >\mu$ as soon as  $m-r<n$ as in this case it is strictly weaker than
\eqref{gaphals}.

\begin{remark}
It is actually possible to solve the equation $F\Psi=\phi$ above with a sharper estimate, by
means of the  complex
\begin{multline*}
\to \Lambda^{r+2}E\otimes S^2Q^*\otimes\det Q^*\to \\ \Lambda^{r+1}E\otimes Q^*\otimes\det Q^*
\to \Lambda^rE\otimes\det Q^*\to\C\to 0
\end{multline*}
with the same mappings as before. Combined with the Fuhrmann trick one can then obtain
Theorem~\ref{bollar}. 
This complex, and the corresponding residues that appear when $f$ is 
only generically surjective, will be studied in a forthcoming paper.
\end{remark}

\def\listing#1#2#3{{\sc #1}:\ {\it #2},\ #3.}

\end{document}